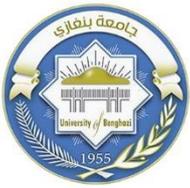

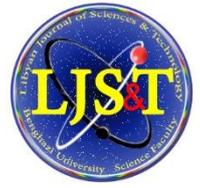

# Existence of at least one positive continuing solution of Urysohn quadratic integral equation by Schauder fixed-point theorem.

**Insaf F. Ben Saoud[a,*], Haitham A. Makhzoum[b], Kheria M. Msaik[c]**

[a]Department of Mathematics, Faculty of Education, University of Benghazi, Benghazi, Libya

[b]Department of Mathematics, Faculty of Science, University of Benghazi, Benghazi, Libya

[c]Department of Mathematics, Faculty of Science, University of Al Zintan, Al Zintan, Libya.

**Highlights**

- **Introduce the concept of the quadratic integral equation and its importance.**
- **We have made several major assumptions that contribute to achieving our goal.**
- **Present and prove the main result.**

A B S T R A C T

We employ Schauder fixed-point Theorem to prove the existence of at least one positive continuous solution of the quadratic integral equation

$$x(t) = a(t) + \int_0^t f_1(t,s,x(s))ds \int_0^t f_2(t,s,x(s))ds, \quad t \in [0,T]$$

Moreover, the maximal and the minimal solutions of the last equation are also proved.

## 1. Introduction

The study of integral equations is one of the most important topics that researchers are interested in, it arises in many scientific fields for instance engineering, and mathematical and scientific analysis. The first who mentioned the term integral equations is Du Bois-Reymond (1888). As a result, a lot of interest appeared from researchers, and the most important of these researchers are Laplace, Fourier, Poission, Liouville, and Able. Upadhyay et al., (2015) provided some special types of integral equations. The quadratic integral equation is a special form of integral equations. The initial study appeared by Chandrasekhar (1947). More appearance of Quadratic integral equation was in the theory of radiative transfer, kinetic theory of gases, in the theory of neutron transport, and in the traffic theory, see Argyros (1985), Banaś et al. (2007), El-Sayed et al. (2008).

Due to the importance of the quadratic integral equation, researchers are interested in studying the existence of its solution, and one of the most important methods that have been used to prove the existence of the solution of the quadratic integral equation is fixed-point theory. Uses of fixed-point theory appeared in many scientific articles for a wealth of reference material on the subject, we refer to Elmabrok et al. (2018), Khalili et al. (2019), Rao, et al. (2020) and the references in them.

Recently, El-Sayed, et al. (2008) considered the quadratic integral equation

$$x(t) = a(t) + \int_0^t f(s,x(s))ds \int_0^t g(s,x(s))ds \qquad (1.1)$$

They proved the existence of at least one continuing positive solution; also, they proved the existence of the maximal and minimal solutions.

Mohamed et al. (2014) discussed the quadratic integral equation

$$x(t) = a(t) + \int_0^t f_1(t,s,x(s))ds \int_0^t f_2(t,s,x(s))ds, \ t \in [0,T] \quad (1.2)$$

and show that it has a unique positive continuous solution by using Banach fixed point Theorem.

The purpose of this study is mainly concerned with at least exists one continuing positive solution of Eq (1.2), and the existence of the solutions of the maximal and minimal by means of the Schauder fixed-point Theorem, in which we provide new conditions that match the requirements of the theory. Also, this article is a scientific addition for every researcher interested in studying quadratic integral equations, as it sheds light on the possibility of using different methods to prove the existence of maximal and minimal solutions to integral equations and obtain the same results according to the standards and conditions imposed for each problem. This article has been organized as follows: Section 2 provides the basic definitions as well as the most important theories of the Time scale which will need it, and Section 3 contains the main results.

## 2. Preliminaries

We present a collection of auxiliary facts in this section which are further required. Let $I = [0, T]$, and $L^1 = L^1[0, T]$ be Lebesgue's Space integrable functions on $I$. Since we are searching for at least exists one continuing positive solution of Eq (1.2), it is natural to assume that

I. $a : I = [0, T] \to R_+$ is continuous, $a = \sup_{t \in [0,T]} |a(t)|$.

II. $f_i : [0, T] \times [0, T] \times R_+ \to R_+$

are Carathéodary functions ( i.e. measurable in $(t, s)$ for all $x \in R_+$ and continuous in $x$ for almost all $(t, s) \in [0, T] \times [0, T]$) and there exist the functions $m_i(t, s)$ such that

$$|f_i(t, s, x)| \leq m_i(t, s), \quad i = 1,2$$

and $\int_0^t m_i(t, s) ds \leq M_i, \quad i = 1,2, \ t \in [0, T]$. Moreover $f_i, \ i = 1,2$ are monotonic nonincreasing in $t \in [0, T]$.

Now, we state the main theories that will have an effective role to reach the desired result

**Theorem 2.1 Schauder fixed-point Theorem (**Geobel *et al.*, 1990**)**

Let $\Psi$ be a convex subset of a Banach space $\mathcal{B}$, Suppose $\mathcal{F} : \Psi \to \Psi$ is compact, continuous. Then $\mathcal{F}$ has at least one fixed-point in $\Psi$.

**Theorem 2.2 Arzela-Ascoli Theorem** (Kolmogorov *et al.*, 1975)

Let $\mho$ be a compact metric space and $C(\mho)$ be the Banach space of real or complex-valued continuous functions normed by

$$\|f\| = \max_{t \in \mho} |f(t)|.$$

If $\aleph = \{f_n\}$ is a sequence in $C(\mho)$ such that $f_n$ is uniformly bounded and equi-continuous. Then the closure of $\aleph$ is compact.

**Theorem 2.3 Lebesgue Dominated Convergence Theorem** (Kolmogorov *et al.*, 1975)

Let $\{\ell_n\}$ be a sequence of functions converging to a limit $\ell$ on $A$, and suppose that

$$|\ell_n(t)| \leq \emptyset(t), t \in A, \ n = 1, 2, 3, \dots$$

Where $\emptyset$ is an integrable function on $A$. Then $\ell$ is integrable on $A$ and

$$\lim_{n \to \infty} \int_A \ell_n(t) \, d\mu = \int_A \ell(t) \, d\mu.$$

Lastly, we provide the next definition, introduced by Lakshmikantham *et al.* (1969) which will be needed later in this paper.

**Definition 2.4** Let $c(t)$ be a solution of the quadratic integral equation Eq (1.1). Then $c(t)$ is said to be a maximal solution of Eq (1.1) if every solution $x(t)$ of Eq (1.1) satisfies the inequality.

$$x(t) < c(t), \quad t \in [0, T]. \tag{2.1}$$

A minimal solution $n(t)$ can be defined similarly by reversing the inequality (2.1) i.e

$$x(t) > n(t), \quad t \in [0, T].$$

## 3. Main result

We introduce and prove the main result in this section. Allow $C = C[0, T]$ to be the continuous functions space on I and set $S$ by

$S = \{x \in C : 0 < x \leq r\} \subset C[0, T]$, where $r = a + M_1 M_2$.

It is obvious that $S$ is closed, convex, bounded, and nonempty.

**Theorem 3.1**

Suppose (I) and (II) are satisfied, then Eq (1.2) has at least one continuing positive solution $x \in C[0, T]$.

**Proof**

Define the mapping $F$ by

$$Fx(t) = a(t) + \int_0^t f_1(t, s, x(s)) ds \int_o^t f_2(t, s, x(s)) ds \tag{3.1}$$

let $x \in S$, then

$$|Fx(t)| = \left| a(t) + \int_0^t f_1(t, s, x(s)) ds \int_o^t f_2(t, s, x(s)) ds \right|$$

$$\leq |a(t)| + \int_0^t |f_1(t, s, x(s))| \, ds \int_0^t |f_2(t, s, x(s))| \, ds$$

$$\leq |a(t)| + \int_0^t m_1(t, s) \, ds \int_0^t m_2(t, s) \, ds$$

$$\leq a + M_1 M_2 = r.$$

This leads that $F$ and $\{F(x)\}$ is uniformly bounded.

Let $t_1, t_2 \in [0, T], t_1 < t_2$ and $|t_2 - t_1| \leq \delta$, then

$$|Fx(t_2) - Fx(t_1)| = |a(t_2) - a(t_1)$$

$$+ \int_0^{t_2} f_1(t_2, s, x(s)) \, ds \int_0^{t_2} f_2(t_2, s, x(s)) ds$$

$$- \int_0^{t_1} f_1(t_1, s, x(s)) \, ds \int_0^{t_1} f_2(t_1, s, x(s)) \, ds \bigg|$$

$$= |a(t_2) - a(t_1)$$

$$+ \left[ \int_0^{t_1} f_1(t_2, s, x(s)) ds + \int_{t_1}^{t_2} f_1(t_2, s, x(s)) ds \right] \int_0^{t_2} f_2(t_2, s, x(s)) ds$$

$$- \int_0^{t_1} f_1(t_1, s, x(s)) \, ds \int_0^{t_1} f_2(t_1, s, x(s)) \, ds \bigg|$$

$$\leq |a(t_2) - a(t_1)|$$

$$+ \left| \int_0^{t_1} f_1(t_2, s, x(s)) ds \int_0^{t_2} f_2(t_2, s, x(s)) ds \right.$$

$$+ \int_{t_1}^{t_2} f_1(t_2, s, x(s)) ds \int_0^{t_2} f_2(t_2, s, x(s)) ds$$

$$- \int_0^{t_1} f_1(t_1, s, x(s)) \, ds \int_0^{t_1} f_2(t_1, s, x(s)) \, ds \bigg|$$

$$\leq |a(t_2) - a(t_1)|$$

$$+ \left| \int_0^{t_1} f_1(t_1, s, x(s)) ds \int_0^{t_2} f_2(t_1, s, x(s)) ds \right.$$

$$+ \int_{t_1}^{t_2} f_1(t_1, s, x(s)) ds \int_0^{t_2} f_2(t_1, s, x(s)) ds$$

$$- \int_0^{t_1} f_1(t_1, s, x(s)) \, ds \int_0^{t_1} f_2(t_1, s, x(s)) \, ds \bigg|$$

$$\leq |a(t_2) - a(t_1)|$$

$$+ \left| \int_0^{t_1} f_1(t_1,s,x(s))ds \left[ \int_0^{t_2} f_2(t_1,s,x(s))ds \right. \right.$$

$$\left. \left. - \int_0^{t_1} f_2(t_1,s,x(s))\,ds \right] \right.$$

$$\left. + \int_{t_1}^{t_2} f_1(t_1,s,x(s))ds \int_0^{t_2} f_2(t_1,s,x(s))ds \right|$$

$$\leq |a(t_2) - a(t_1)|$$

$$+ \int_0^{t_1}|f_1(t_1,s,x(s))|ds \int_{t_1}^{t_2}|f_2(t_1,s,x(s))|ds$$

$$+ \int_0^{t_2}|f_2(t_1,s,x(s))|ds \int_{t_1}^{t_2}|f_1(t_1,s,x(s))|ds$$

$$\leq |a(t_2) - a(t_1)|$$

$$+ \int_0^{t_1} m_1(t,s)\,ds \int_{t_1}^{t_2} m_2(t,s)\,ds$$

$$+ \int_0^{t_2} m_2(t,s)\,ds \int_{t_1}^{t_2} m_1(t,s)\,ds.$$

That tells us $F\{x\}$ is eque-continuous on $[0,T]$. Using Arzela-Ascoli Theorem, we found that $F$ is compact.

Now we show that $F: S \to S$ is continuous. Let $\{x_n\} \subset S$, and $x_n \to x$, then

$$Fx_n(t) = a(t) + \int_0^t f_1(t,s,x_n(s))\,ds \int_0^t f_2(t,s,x_n(s))\,ds$$

$$\lim_{n \to \infty} Fx_n(t) = \lim_{n \to \infty} a(t) + \lim_{n \to \infty}\{\int_0^t f_1(t,s,x_n(s))ds \int_0^t f_2(t,s,x_n(s))ds\}.$$

Now

$$f_i(t,s,x_{n_k}) \to f_i(t,s,x), \quad i=1,2.$$

Also

$$|f_i(t,s,x_{n_k})| \leq m_i(t,s), \quad i=1,2.$$

By the use of Theorem 2.3, we have

$$Fx(t) = \lim_{n_k \to \infty} Fx_{n_k}(t)$$

$$= a(t) + \int_0^t \lim_{n_k \to \infty} f_1(t,s,x_{n_k}(s))\,ds \int_0^t \lim_{n_k \to \infty} f_2(t,s,x_{n_k}(s))\,ds$$

$$Fx(t) = a(t) + \int_0^t f_1(t,s,x(s))ds \int_0^t f_2(t,s,x(s))ds.$$

Then $Fx_n(t) \to Fx(t)$. Which leads to $F$ be continuous.

Since all requisites of Theorem 2.1 are satisfied, then the mapping $F$ has at least one continuing positive solution $x \in C[0,T]$, the proof is complete.

**Corollary 3.2** Assume that $f,g:[0,T] \times R_+ \to R_+$ are $L^1$-Carathéodory functions with $|f| \leq m_1$ and $|g| \leq m_2$, then Eq (1.1) has at least one continuing positive solution.

Now we want to prove the existence of maximal and minimal solutions. We will present the following lemma.

**Lemma 3.3** Let $f_i(t,s,x)$, $i=1,2$ satisfy the assumption (II) and $x(t), y(t)$ are two continuous functions on $[0,T]$ satisfying

$$x(t) \leq a(t) + \int_0^t f_1(t,s,x(s))ds \int_0^t f_2(t,s,x(s))ds \quad, t \in [0,T]$$

$$y(t) \geq a(t) + \int_0^t f_1(t,s,y(s))ds \int_0^t f_2(t,s,y(s))ds \quad, t \in [0,T]$$

and one of them is strict.

If $f_i$, $i=1,2$ are monotonic non-decreasing in $x$, then

$$x(t) < y(t), \quad t > 0 \tag{3.2}$$

Let the conclusion (3.2) be false, then there exists $t_1$ such that

$$x(t_1) = y(t_1), \quad t_1 > 0$$

and

$$x(t) < y(t), \quad 0 < t < t_1.$$

From the monotonicity of $f_1, f_2$ in $x$, we get

$$x(t_1) \leq a(t_1) + \int_0^{t_1} f_1(t_1,s,x(s))ds \int_0^{t_1} f_2(t_1,s,x(s))ds, \quad t \in [0,T]$$

$$< a(t_1) + \int_0^{t_1} f_1(t_1,s,y(s))ds \int_0^{t_1} f_2(t_1,s,y(s))ds, \quad t \in [0,T]$$

$x(t_1) < y(t_1)$ which contradicts the fact that $x(t_1) = y(t_1)$. Then

$$x(t) < y(t).$$

Now, we will introduce the next theorem to show the maximal and minimal solutions of Eq (1.2).

**Theorem 3.4**

Let the assumptions (I) and (II) of Theorem 2.1 be satisfied. If $f_1(t,s,x), f_2(t,s,x)$ are monotonic non-decreasing in $x$ for each $t \in [0,T]$, then the maximal and minimal solutions of Eq (1.2) exist.

**Proof.**

Firstly, we will prove the existence of the maximal solution of Eq (1.2). Let $\epsilon > 0$ be given, and consider

$$x_\epsilon(t) \leq a(t) + \int_0^t f_{1_\epsilon}(t,s,x_\epsilon(s))ds \int_0^t f_{2_\epsilon}(t,s,x_\epsilon(s))ds, \quad t \in [0,T] \tag{3.3}$$

where

$$f_{i_\epsilon}(t,s,x_\epsilon(t)) = f_i(t,s,x_\epsilon(t)) + \epsilon, \quad i=1,2.$$

Clearly, the functions $f_{i_\epsilon}(t,s,x_\epsilon(t))$, $i=1,2$ are $L^1$- Carathéodory functions, therefore (3.3) has a solution on $C[0,T]$. Let $\epsilon_1, \epsilon_2$ such that $0 < \epsilon_2 < \epsilon_1 < \epsilon$, then

$$x_{\epsilon_2}(t) = a(t) + \int_0^t f_{1_{\epsilon_2}}(t,s,x_{\epsilon_2}(s))\,ds \int_0^t f_{2_{\epsilon_2}}(t,s,x_{\epsilon_2}(s))\,ds$$

$$= a(t) + \int_0^t (f_1(t,s,x_{\epsilon_2}(s)) + \epsilon_2)ds \int_0^t (f_2(t,s,x_{\epsilon_2}(s)) + \epsilon_2)ds \tag{3.4}$$

also

$$x_{\epsilon_1}(t) = a(t) + \int_0^t f_{1_{\epsilon_1}}(t,s,x_{\epsilon_1}(s))\,ds \int_0^t f_{2_{\epsilon_1}}(t,s,x_{\epsilon_1}(s))\,ds$$

$$= a(t) + \int_0^t (f_1(t,s,x_{\epsilon_1}(s)) + \epsilon_1)ds \int_0^t (f_2(t,s,x_{\epsilon_1}(s)) + \epsilon_1)ds$$

$$x_{\epsilon_1} > a(t) + \int_0^t f_1\big(t,s,x_{\epsilon_1}(s)\big)ds \int_0^t f_2\big(t,s,x_{\epsilon_1}(s)\big)ds \quad (3.5)$$

Applying Lemma 3.3 on (3.4) and (3.5), we have

$$x_{\epsilon_2}(t) < x_{\epsilon_1}(t) \quad for\ t \in [0,T].$$

As previously shown, $x_\epsilon(t)$ is equi-continuous and uniformly bounded. Then, by using Arzela- Ascoli theorem, a decreasing sequence $\epsilon_n$ exists such that $\epsilon_n \to 0$ as $n \to \infty$, and $\lim_{n\to\infty} x_{\epsilon_n}(t)$ exists uniformly in $[0,T]$ and denote the limit by $q(t)$.

Since $f_{i_\epsilon}\big(t,s,x_\epsilon(t)\big)$, $i = 1,2$ is continued in the third argument, we get

$$f_{i_\epsilon}\big(t,s,x_\epsilon(t)\big) \to f_i(t,s,x(t)) \text{ as } n \to \infty, \quad i = 1,2.$$

And

$$q(t) = \lim_{n\to\infty} x_{\epsilon_n}(t) = a(t) + \int_0^t f_1(t,s,q(s))ds \int_0^t f_2(t,s,q(s))ds$$

which implies that $q(t)$ is a solution of Eq (1.2).

Finally, we will prove that $q(t)$ has the maximal solution of Eq (1.2). let $x(t)$ be any solution of Eq (1.2), then

$$x(t) = a(t) + \int_0^t f_1(t,s,x(s))ds \int_0^t f_2(t,s,x(s))ds \quad (3.6)$$

also

$$x_\epsilon(t) = a(t) + \int_0^t f_{1_\epsilon}(t,s,x_\epsilon(s))ds \int_0^t f_{2_\epsilon}(t,s,x_\epsilon(s))ds$$

$$x_\epsilon(t) = a(t) + \int_0^t (f_1(t,s,x_\epsilon(s)) + \epsilon)ds \int_0^t (f_2(t,s,x_\epsilon(s)) + \epsilon)ds$$

$$x_\epsilon(t) > a(t) + \int_0^t f_1(t,s,x_\epsilon(s))ds \int_0^t f_2(t,s,x_\epsilon(s))ds. \quad (3.7)$$

Applying Lemma 3.3 on (3.6) and (3.7), we get

$$x(t) < x_\epsilon(t), \quad for\ t \in [0,T].$$

it is clear that $x_\epsilon(t)$ tends to q(t) uniformly in [0, T] as $\epsilon \to \infty$, from the uniqueness of the maximal solution, see Lakshmikantham, *et al.* (1969). By the same steps, we can show the Eq (1.2) has the minimal solution. We set that

$$f_{i_\epsilon}\big(t,s,x_\epsilon(t)\big) = f_i(t,s,x_\epsilon(t)) - \epsilon, \quad i = 1,2$$

and show the minimal solution exists.

### Corollary 3.5

Let the functions $f$ and $g$ be non-decreasing in the second argument and the assumptions of Corollary 3.2 are satisfied. Then Eq (1.1) has maximal and minimal solutions.

## 4. Conclusion

This article contains three main parts. The first part is a historical overview of the topic and shows the importance of the quadratic integral equations. The second part details the assumptions and theories used to obtain the desired result. The final part provides the main result. As introducing Theorem 3.1 to attain the theory's conditions, we use Theorem 2.2 that showed that $F$ in Eq (3.1) is compact. We also employ Theorem 2.3 to prove that $F$ is continuous, which helps us achieve at least one continuing positive solution of Eq (1.2). Furthermore, we provide Theorem 3.4 to show the solutions of the maximal and minimal of Eq (1.2) with the help of Lemma 3.3 and Theorem 2.2.

### Acknowledgment

The authors express their appreciation to the referees for reading, corrections, and valuable hints of this manuscript.